# On Closed Space Curves in Mınkowskı Space-Time $E_v^n$

**Zehra ARI, Mehmet ÖNDER**
*Celal Bayar University, Department of Mathematics, Faculty of Arts and Sciences, Manisa, Turkey.*
E-mails: zehra.ari@bayar.edu.tr, mehmet.onder@bayar.edu.tr

**Abstract**
In this paper, we give the generalization of the criterion for a 3-space curve to be closed given by [3] to an $n$-space curve in Minkowski space-time $E_v^n$. Furthermore, we apply this criterion for a curve lying on an oriented surface in the Minkowski space $E_1^3$ as an application.

**Key words:** Minkowski space-time; closed curves; Serre-Frenet frame.
**AMS Subject Classification:** 31A35, 53C50, 53B30.

## 1. Introduction

In 1978, Professor Shiing-Shen Chern raised the following question at the Institute of Mathematics of Academia Sinica:
*"What is the necessary and sufficient condition to be satisfied by the curvature and torsion, so that spaces curve to be a closed one?"*
The answer of the question is given by H. C. Chung in 1981. By considering the Frenet frame of a space curve in Euclidean 3-space $E^3$ and using the well-known method of successive approximation due to E. Picard, Chung presented a criterion for a 3-space curve to be closed[1]. By using this criterion, he found that the curve determined by curvature $\kappa(s) = cons. > 0$ and torsion $\sigma(s) = const.$ is periodic of period $\omega$, if and only if $\sigma = 0$ and $\kappa = 2\pi/\omega$. Köse and et al generalized the method of Chung to the space curves in $n$-dimensional Euclidean space $E^n$ [7].

In [1], Altin gave the Serre-Frenet frame and harmonic curvatures of non-null curves in Minkowski space-time $E_v^n$. Furthermore, the geometry of the curves lying on an oriented surface has been given in Minkowski 3-space $E_1^3$ in [6,11,12,13].

In this paper, we give a criterion for closed space curves in Minkowski space-time $E_v^n$. In addition to this, as a special case we considered a curve lying on a surface in Minkowski 3-space $E_1^3$ and by regarding the Darboux frame of this curve we obtained a necessary and sufficient condition to be satisfied by the geodesic curvature, normal curvature and geodesic torsion.

## 2. Preliminaries

Minkowski space-time $E_v^n$ is an Euclidean space $E^n$ provided with the standard flat metric given by
$$\langle,\rangle = -\sum_{i=1}^{v} dx_i^2 + \sum_{i=v+1}^{n} dx_i^2$$
where $(x_1, x_2, ..., x_n)$ is a rectangular coordinate system in $E_v^n$. Since $\langle,\rangle$ is an indefinite metric, recall that a vector $\vec{v} \in E_v^n$ can have one of three causal characters; it can be spacelike if $\langle \vec{v}, \vec{v} \rangle > 0$ or $\vec{v} = 0$, timelike if $\langle \vec{v}, \vec{v} \rangle < 0$ and null(lightlike) if $\langle \vec{v}, \vec{v} \rangle = 0$ and $\vec{v} \neq 0$. Let $x(s): I \subset \mathbb{R} \to E_v^n$ be a regular curve in $E_v^n$. Analogues to the vectors, the curve $x(s)$ in $E_v^n$ can locally be spacelike, timelike or null (lightlike), if all of its velocity vectors $x'(s)$ are spacelike, timelike or null (lightlike), respectively. A timelike (resp. spacelike) curve $x(s)$ is said to be parameterized by a pseudo-arclength parameter $s$, if $\langle x'(s), x'(s) \rangle = -1$ (resp.

$\langle x'(s), x'(s) \rangle = 1$). Also, recall that the pseudo-norm of an arbitrary vector $\vec{v} \in E_v^n$ is given by $\|\vec{v}\| = \sqrt{|\langle \vec{v}, \vec{v} \rangle|}$. Therefore $\vec{v}$ is a unit vector if $\langle \vec{v}, \vec{v} \rangle = \pm 1$. The velocity of the curve $x(s)$ is given by $\|x'(s)\|$. Next, vectors $\vec{v}$, $\vec{w}$ in $E_v^n$ are said to be orthogonal if $\langle \vec{v}, \vec{w} \rangle = 0$.

Let $(\vec{e}_1, \vec{e}_2, ..., \vec{e}_n)$ denote an orthonormal basis in the Minkowski space-time $E_v^n$. Then the position vector of a curve $x(s)$ can be given by $\vec{x}(s) = x_1\vec{e}_1 + x_2\vec{e}_2 + \cdots + x_n\vec{e}_n$. Let $\omega$ be the total length of the curve $x(s)$. When $x(s)$ is a closed curve, the position vector $\vec{x}(s)$ must be a periodic function of period $\omega$[7].

Let $x(s)$ be a regular curve in $E_v^n$ and $\varphi = \{x'(s), x''(s), ..., x^{(n)}(s)\}$ a maximal linear independent and non-null set. The orthonormal system $\{\vec{V}_1(s), \vec{V}_2(s), ..., \vec{V}_n(s)\}$ can be obtained from $\varphi$. This system is called a moving Serre-Frenet frame along the curve $x(s)$ in the space $E_v^n$. In this paper, we consider the curve $x(s)$ whose derivatives $x^{(l)}(s)$, $(1 \leq l \leq n)$ are non-null.

**Definition 2.1.** Let $x(s)$ be a regular curve in $E_v^n$ and $\{\vec{V}_1(s), \vec{V}_2(s), ..., \vec{V}_n(s)\}$ denote the Frenet frame of $x(s)$. The functions $k_i : I \to \mathbb{R}$ defined by

$$k_i(s) = \varepsilon_i(s)\varepsilon_{i-1}(s)\langle \vec{V}_i'(s), \vec{V}_{i+1}'(s) \rangle, \quad 1 \leq i \leq n-1,$$

are called curvature functions of $x(s)$. Here $\varepsilon_i(s) = \langle \vec{V}_i(s), \vec{V}_i(s) \rangle = \pm 1$. Furthermore, the real number $k_i(s)$ is called the $i$-th curvature on $x$ at the point $x(s)$ [1].

**Theorem 2.1.** *Let $x(s)$ be a unit speed curve in $E_v^n$ and let the set $\{\vec{V}_1(s), \vec{V}_2(s), ..., \vec{V}_n(s)\}$ denote the Serre-Frenet frame at the point $x(s)$. Then, the derivative formulae are given as follows*

$$\vec{V}_1' = \varepsilon_1(s)k_1(s)\vec{V}_2$$
$$\vec{V}_i' = -\varepsilon_i(s)k_{i-1}(s)\vec{V}_{i-1}(s) + \varepsilon_i(s)k_i(s)\vec{V}_{i+1}(s), \quad (1 < i < n)$$
$$\vec{V}_n' = -\varepsilon_r(s)k_{n-1}(s)\vec{V}_{n-1}(s)$$

*or in the matrix form*

$$\begin{bmatrix} V_1' \\ V_2' \\ V_3' \\ \vdots \\ V_{n-1}' \\ V_n' \end{bmatrix} = \begin{bmatrix} 0 & \varepsilon_1 k_1 & 0 & 0 & \cdots & 0 & 0 & 0 \\ -\varepsilon_2 k_1 & 0 & \varepsilon_2 k_2 & 0 & \cdots & 0 & 0 & 0 \\ 0 & -\varepsilon_3 k_2 & 0 & \varepsilon_3 k_3 & \cdots & 0 & 0 & 0 \\ \vdots & \vdots & \vdots & \vdots & \vdots & \vdots & \vdots & \vdots \\ 0 & 0 & 0 & 0 & \cdots & -\varepsilon_{n-1} k_{n-2} & 0 & \varepsilon_{n-1} k_{n-1} \\ 0 & 0 & 0 & 0 & \cdots & 0 & -\varepsilon_n k_{n-1} & 0 \end{bmatrix} \begin{bmatrix} V_1 \\ V_2 \\ V_3 \\ \vdots \\ V_{n-1} \\ V_n \end{bmatrix}. \quad (1)$$

[1]. In the special case $v = 1$, $n = 3$ we have Minkowski 3-space $E_1^3$ and we can give the followings for this space.

**Definition 2.2.** A surface in the Minkowski 3-space $E_1^3$ is called a timelike surface if the induced metric on the surface is a Lorentz metric and is called a spacelike surface if the induced metric on the surface is a positive definite Riemannian metric, i.e., the normal vector on the spacelike (timelike) surface is a timelike (spacelike) vector[5].

Let $S$ be an oriented surface in three-dimensional Minkowski space $E_1^3$ and let consider a non-null curve $x(s)$ lying on $S$ fully. Since the curve $x(s)$ lies on the surface $S$ there exists another frame of the curve $x(s)$ which is called *Darboux frame* and denoted by $\{\vec{T}, \vec{g}, \vec{n}\}$. In this frame $\vec{T}$ is the unit tangent of the curve, $\vec{n}$ is the unit normal of the surface $S$ and $\vec{g}$ is a unit vector given by $\vec{g} = \vec{n} \times \vec{T}$.

According to the Lorentzian causal characters of the surface $S$ and the curve $x(s)$ lying on $S$, the derivative formulae of the Darboux frame can be changed as follows:

**i)** If the surface $S$ is a timelike surface, then the curve $x(s)$ lying on $S$ can be a spacelike or a timelike curve. Thus, the derivative formulae of the Darboux frame of $x(s)$ is given by

$$\begin{bmatrix} \dot{\vec{T}} \\ \dot{\vec{g}} \\ \dot{\vec{n}} \end{bmatrix} = \begin{bmatrix} 0 & k_g & -\varepsilon k_n \\ k_g & 0 & \varepsilon \tau_g \\ k_n & \tau_g & 0 \end{bmatrix} \begin{bmatrix} \vec{T} \\ \vec{g} \\ \vec{n} \end{bmatrix}, \quad \langle \vec{T}, \vec{T} \rangle = \varepsilon = \pm 1, \ \langle \vec{g}, \vec{g} \rangle = -\varepsilon, \ \langle \vec{n}, \vec{n} \rangle = 1. \tag{2}$$

**ii)** If the surface $S$ is a spacelike surface, then the curve $x(s)$ lying on $S$ is a spacelike curve. Thus, the derivative formulae of the Darboux frame of $x(s)$ is given by

$$\begin{bmatrix} \dot{\vec{T}} \\ \dot{\vec{g}} \\ \dot{\vec{n}} \end{bmatrix} = \begin{bmatrix} 0 & k_g & k_n \\ -k_g & 0 & \tau_g \\ k_n & \tau_g & 0 \end{bmatrix} \begin{bmatrix} \vec{T} \\ \vec{g} \\ \vec{n} \end{bmatrix}, \quad \langle \vec{T}, \vec{T} \rangle = 1, \ \langle \vec{g}, \vec{g} \rangle = 1, \ \langle \vec{n}, \vec{n} \rangle = -1. \tag{3}$$

In these formulae $k_g$, $k_n$ and $\tau_g$ are called the geodesic curvature, the normal curvature and the geodesic torsion, respectively [6,11,12,13].

## 3. Closed Space Curves in Minkowski space-time $E_v^n$

In this section, by considering the Serre-Frenet frame of the curve, we give a criterion for closed space curves in Minkowski space-time $E_v^n$. For this purpose, we use the method given in [3,7].

Firstly, let consider the system of linear differential equations

$$\frac{d\phi_i}{dt} = \sum_{j=1}^{n} a_{ij}(s)\phi_j, \quad (i = 1, 2, ..., n), \tag{4}$$

where $a_{ij}(s)$ are assumed to be continuous and periodic of period $\omega$. Let the initial conditions be $\phi_i(0) = \lambda_i$. The matrix form of (4) can be given by

$$\frac{d\phi}{dt} = A(s)\phi, \quad \phi = \begin{bmatrix} \phi_1 \\ \phi_2 \\ \vdots \\ \phi_n \end{bmatrix}, \quad A(s) = \begin{bmatrix} a_{11} & a_{12} & \cdots & a_{1n} \\ a_{21} & a_{22} & \cdots & a_{2n} \\ \cdots & \cdots & \cdots & \cdots \\ a_{n1} & a_{n2} & \cdots & a_{nn} \end{bmatrix}.$$

Since the $a_{ij}(s)$ are periodic, $A(s)$ is a continuous periodic matrix function of period $\omega$. Let us write $\int_0^s A(s)ds = \left( \int_0^s a_{ij}(s)ds \right)$ an $n \times n$ matrix and

$$\lambda = \begin{bmatrix} \lambda_1 \\ \lambda_2 \\ \vdots \\ \lambda_n \end{bmatrix}, \quad \phi(s,\lambda) = \begin{bmatrix} \phi_1(s,\lambda) \\ \phi_2(s,\lambda) \\ \vdots \\ \phi_n(s,\lambda) \end{bmatrix}.$$

Then the equation (4) may be abbreviated to the form

$$\frac{d\phi}{ds} = A(s)\phi, \quad \phi(0) = \lambda. \tag{5}$$

[See 4,9]. From (5) we get the integral equation $\phi(s) = \lambda + \int_0^s A(\sigma)\phi(\sigma)d\sigma$ [See 4,10]. The solution $\phi(s,\lambda)$ of the equation $\frac{d\phi}{ds} = A(s)\phi, \ \phi(0) = \lambda, \ \big(A(s+\omega) \equiv A(s)\big)$ is periodic if and only if $\int_0^\omega A(\sigma)\phi(\sigma,\lambda)d\sigma = 0$ which is given by a lemma in [3] as follows:

**Lemma 2.1.** *The solution $\phi(s,\lambda)$ of the equation $\frac{d\phi}{ds} = A(s)\phi, \ \phi(0) = \lambda, \ \big(A(s+\omega) \equiv A(s)\big)$ is periodic of period $\omega$ if and only if $\int_0^\omega A(\sigma)\phi(\sigma,\lambda)d\sigma = 0$.*

In the same paper, using the well-known method of successive approximation due to E. Picard, the solution $\phi(s,\lambda)$ is constructed as

$$\phi(s,\lambda) = \left\{ I + \xi A(s) + \xi^{(2)} A(s) + \cdots + \xi^{(n)} A(s) + \cdots \right\} \lambda$$

where $I$ is $n \times n$ unit matrix. The following theorem is proved by applying Lemma 2.1 to the expression

$$\int_0^\omega A(s)\phi(s,\lambda)ds = \left\{ \xi A(s) + \xi^{(2)} A(s) + \cdots + \xi^{(n)} A(s) + \cdots \right\} \lambda = M(s)\lambda$$

where

$$\xi A(s) = \xi^{(1)} A(s) = \int_0^s A(\sigma)d\sigma, \quad \xi^{(n)} A(s) = \int_0^s A(\sigma)\xi^{(n-1)} A(\sigma)d\sigma, \ (n>1)$$

and

$$M(s) = \xi A(s) + \xi^{(2)} A(s) + \cdots + \xi^{(n)} A(s) + \cdots$$

**Theorem 3.1.** *The equations $\frac{d\phi}{ds} = A(s)\phi$ have a non-vanishing periodic solution of period $\omega$ if and only if $\det(M(\omega)) = 0$. In particular, the equations $\frac{d\phi}{ds} = A(s)\phi$ have $n$-linearly independent periodic solutions of period $\omega$ if and only if the matrix $M(\omega) = 0$*[3].

Let now consider the Serre-Frenet formulae given in equations (1). If we write the Serre-Frenet vectors $\vec{V}_1, \vec{V}_2, \ldots \vec{V}_n$ in coordinates form we get

$$\vec{V}_i = \sum_{j=1}^n v_{ij}\vec{e}_j, \ (i = 1,2,\ldots,n). \tag{6}$$

From (1) and (6) we obtain the systems of linear differential equations:

$$\begin{cases} \dfrac{dv_{1j}}{ds} = \varepsilon_1 k_1(s) v_{2j} \\ \dfrac{dv_{2j}}{ds} = -\varepsilon_2 k_1(s) v_{1j} + \varepsilon_2 k_2(s) v_{3j} \\ \vdots \\ \dfrac{dv_{(n-1)j}}{ds} = -\varepsilon_{n-1} k_{n-2}(s) v_{(n-2)j} + \varepsilon_{n-1} k_{n-1}(s) v_{nj} \\ \dfrac{dv_{nj}}{ds} = -\varepsilon_n k_{n-1}(s) v_{(n-1)j}, \quad (j=1,2,\ldots,n) \end{cases} \quad (7)$$

Thus, we observe that

$(v_{11}, v_{21}, \ldots, v_{n1})$
$(v_{12}, v_{22}, \ldots, v_{n2})$
$\vdots$
$(v_{1n}, v_{2n}, \ldots, v_{nn})$

are $n$-independent solutions of the following system of differential equations

$$\begin{cases} \dfrac{d\phi_1}{ds} = \varepsilon_1 k_1(s) \phi_2 \\ \dfrac{d\phi_2}{ds} = -\varepsilon_2 k_1(s) \phi_1 + \varepsilon_2 k_2(s) \phi_3 \\ \vdots \\ \dfrac{d\phi_{n-1}}{ds} = -\varepsilon_{n-1} k_{n-2}(s) \phi_{n-2} + \varepsilon_{n-1} k_{n-1}(s) \phi_n \\ \dfrac{d\phi_n}{ds} = -\varepsilon_n k_{n-1}(s) \phi_{n-1} \end{cases} \quad (8)$$

The system in (8) is a special case of the general equations given in (4).

Now, we can give the following theorem which gives a criterion for space curve in Minkowski space-time $E_v^n$ to be closed.

**Theorem 3.2.** Let $x(s)$ be a curve in Minkowski space-time $E_v^n$ with curvatures $k_1(s), k_2(s), \ldots, k_{n-1}(s)$ and let the curvatures $k_1(s), k_2(s), \ldots, k_{n-1}(s)$ be continuous periodic functions of period $\omega$. Then $x(s)$ is a periodic curve with period $\omega$ if and only if

i) the matrix $M(\omega) = 0$,

ii) $\omega + \int_0^\omega m_{11} ds = \int_0^\omega m_{12} ds = \int_0^\omega m_{13} ds = \ldots = \int_0^\omega m_{1n} ds = 0$

where

$$M(t) = \xi A(t) + \xi^{(2)} A(t) + \cdots + \xi^{(n)} A(t) + \cdots,$$

$$A(t) = \begin{bmatrix} 0 & \varepsilon_1 k_1(t) & 0 & 0 & \cdots & 0 & 0 & 0 \\ -\varepsilon_2 k_1(t) & 0 & \varepsilon_2 k_2(t) & 0 & \cdots & 0 & 0 & 0 \\ 0 & -\varepsilon_3 k_2(t) & 0 & \varepsilon_3 k_3(t) & \cdots & 0 & 0 & 0 \\ \vdots & \vdots & \vdots & \vdots & \vdots & \vdots & \vdots & \vdots \\ 0 & 0 & 0 & 0 & \cdots & -\varepsilon_{n-1} k_{n-2}(t) & 0 & \varepsilon_{n-1} k_{n-1}(t) \\ 0 & 0 & 0 & 0 & \cdots & 0 & -\varepsilon_n k_{n-1}(t) & 0 \end{bmatrix}$$

and $m_{ij}(s)$ are the entries of the matrix $M(t)$.

**Proof:** By considering Theorem 3.1 for the system (8) we have

$$A(s) = \begin{bmatrix} 0 & \varepsilon_1 k_1(s) & 0 & 0 & \cdots & 0 & 0 & 0 \\ -\varepsilon_2 k_1(s) & 0 & \varepsilon_2 k_2(s) & 0 & \cdots & 0 & 0 & 0 \\ 0 & -\varepsilon_3 k_2(s) & 0 & \varepsilon_3 k_3(s) & \cdots & 0 & 0 & 0 \\ \vdots & \vdots & \vdots & \vdots & \vdots & \vdots & \vdots & \vdots \\ 0 & 0 & 0 & 0 & \cdots & -\varepsilon_{n-1} k_{n-2}(s) & 0 & \varepsilon_{n-1} k_{n-1}(s) \\ 0 & 0 & 0 & 0 & \cdots & 0 & -\varepsilon_n k_{n-1}(s) & 0 \end{bmatrix}$$

$$\xi A(s) = \begin{bmatrix} 0 & \int_0^s \varepsilon_1 k_1(t) dt & 0 & 0 & \cdots & 0 & 0 & 0 \\ -\int_0^s \varepsilon_2 k_1(t) dt & 0 & \int_0^s \varepsilon_2 k_2(t) dt & 0 & \cdots & 0 & 0 & 0 \\ 0 & -\int_0^s \varepsilon_3 k_2(t) dt & 0 & \int_0^s \varepsilon_3 k_3(t) dt & \cdots & 0 & 0 & 0 \\ \vdots & \vdots & \vdots & \vdots & \vdots & \vdots & \vdots & \vdots \\ 0 & 0 & 0 & 0 & \cdots & -\int_0^s \varepsilon_{n-1} k_{n-2}(t) dt & 0 & \int_0^s \varepsilon_{n-1} k_{n-1}(t) dt \\ 0 & 0 & 0 & 0 & \cdots & 0 & -\int_0^s \varepsilon_n k_{n-1}(t) dt & 0 \end{bmatrix}$$

and so on. When and only when the matrix $M(\omega)$ is a zero matrix, there exist $n$-orthonormal vector functions $\vec{V}_1(s), \vec{V}_2(s), ..., \vec{V}_n(s)$ of period $\omega$ such that each set of functions $\{v_{1j}, v_{2j}, ..., v_{nj}\}$, $(j=1,2,...,n)$ forms a solution of the equation $\dfrac{d\phi}{ds} = A(s)\phi$ corresponding to the initial condition $(a_{1j}, a_{2j}, ..., a_{nj})$.

The vector function defining the curve $\vec{x}(s) = \int_0^\omega \vec{V}_1(s) ds$ where $\vec{V}_1(s)$ is given by

$$\begin{bmatrix} v_{1j} \\ v_{2j} \\ \vdots \\ v_{nj} \end{bmatrix} = (I + M(s)) \begin{bmatrix} a_{1j} \\ a_{2j} \\ \vdots \\ a_{nj} \end{bmatrix}, \quad (j=1,2,...,n)$$

and $I$ is $n \times n$ unit matrix. The curve $x(s)$ is periodic of period $\omega$ if and only if $\int_0^\omega \vec{V}_1(s) ds = 0$. Let now the initial condition be

$$v_{1j}(0) = a_{1j}, \; v_{2j}(0) = a_{2j}, ..., v_{nj}(0) = a_{nj}, \; (j=1,2,...,n)$$

where $(a_{11}, a_{12}, ..., a_{1n})$, $(a_{21}, a_{22}, ..., a_{2n})$, ..., $(a_{n1}, a_{n2}, ..., a_{nn})$ form an orthonormal frame. Then

$$\begin{bmatrix} v_{1j} \\ v_{2j} \\ \vdots \\ v_{nj} \end{bmatrix} = (I + M(t)) \begin{bmatrix} a_{1j} \\ a_{2j} \\ \vdots \\ a_{nj} \end{bmatrix}, \quad (j=1,2,...,n).$$

That is, $v_{1j} = a_{1j} + a_{1j}m_{11} + a_{2j}m_{12} + \cdots + a_{nj}m_{1j}$. So that

$$\int_0^\omega v_{1j}(s)ds = a_{1j} + a_{1j}\int_0^\omega m_{11}ds + \cdots + a_{nj}\int_0^\omega m_{1j}(s)ds.$$

Since the determinant

$$\begin{vmatrix} a_{11} & a_{21} & \cdots & a_{n1} \\ a_{12} & a_{22} & \cdots & a_{n2} \\ \vdots & \vdots & \vdots & \vdots \\ a_{1n} & a_{2n} & \cdots & a_{nn} \end{vmatrix} \neq 0.$$

The condition $\int_0^\omega \vec{V}_1(s)ds = 0$ is equivalent to

$$\omega + \int_0^\omega m_{11}ds = \int_0^\omega m_{12}ds = \int_0^\omega m_{13}ds = \ldots = \int_0^\omega m_{1n}ds = 0.$$

This completes the proof.

In the following example, this criterion is applied to a curve lying on a surface in Minkowski 3-space $E_1^3$ as a special case.

**Example 3.1:** Let $S$ be an oriented timelike surface in three-dimensional Minkowski space $E_1^3$ and let consider a non-null curve $x(s)$ lying on $S$ fully with Darboux frame $\{\vec{T}, \vec{g}, \vec{n}\}$. Assume that $k_g = const.$, $k_n = const.$ and $\tau_g = const.$ for $x(s)$. By considering (2), from Theorem 3.2, we get

$$m_{11}(\omega) = (k_g^2 - \varepsilon k_n^2)\left[\frac{\omega^2}{2!} + (k_g^2 - \varepsilon k_n^2 + \varepsilon \tau_g^2)\frac{\omega^4}{4!} + (k_g^2 - \varepsilon k_n^2 + \varepsilon \tau_g^2)^2 \frac{\omega^6}{6!} + \cdots\right]$$

$$= (k_g^2 - \varepsilon k_n^2)\frac{\cosh\left((k_g^2 - \varepsilon k_n^2 + \varepsilon \tau_g^2)^{1/2}\omega\right) - 1}{k_g^2 - \varepsilon k_n^2 + \varepsilon \tau_g^2}$$

$$m_{12}(\omega) = k_g\left[\omega + (k_g^2 - \varepsilon k_n^2 + \varepsilon \tau_g^2)\frac{\omega^3}{3!} + (k_g^2 - \varepsilon k_n^2 + \varepsilon \tau_g^2)^2 \frac{\omega^5}{5!} + \cdots\right]$$

$$- k_n \tau_g \left[\varepsilon \frac{\omega^2}{2!} - \varepsilon(k_g^2 - \varepsilon k_n^2 + \varepsilon \tau_g^2)\frac{\omega^4}{4!} + (k_g^2 - \varepsilon k_n^2 + \varepsilon \tau_g^2)^2 \frac{\omega^6}{6!} + \cdots\right]$$

$$= k_g \frac{\sinh\left((k_g^2 - \varepsilon k_n^2 + \varepsilon \tau_g^2)^{1/2}\right)\omega}{(k_g^2 - \varepsilon k_n^2 + \varepsilon \tau_g^2)^{1/2}} - \varepsilon k_n \tau_g \frac{\cosh\left((k_g^2 - \varepsilon k_n^2 + \varepsilon \tau_g^2)^{1/2}\omega\right) - 1}{k_g^2 - \varepsilon k_n^2 + \varepsilon \tau_g^2},$$

$$m_{13}(\omega) = -\varepsilon k_n\left[\omega + (k_g^2 - \varepsilon k_n^2 + \varepsilon \tau_g^2)\frac{\omega^3}{3!} + (k_g^2 - \varepsilon k_n^2 + \varepsilon \tau_g^2)^2 \frac{\omega^5}{5!} + \cdots\right]$$

$$+ \varepsilon k_g \tau_g \left[\frac{\omega^2}{2!} + (k_g^2 - \varepsilon k_n^2 + \varepsilon \tau_g^2)\frac{\omega^4}{4!} + (k_g^2 - \varepsilon k_n^2 + \varepsilon \tau_g^2)^2 \frac{\omega^6}{6!} + \cdots\right]$$

$$= -\varepsilon k_n \frac{\sinh\left((k_g^2 - \varepsilon k_n^2 + \varepsilon \tau_g^2)^{1/2}\right)\omega}{(k_g^2 - \varepsilon k_n^2 + \varepsilon \tau_g^2)^{1/2}} + \varepsilon k_g \tau_g \frac{\cosh\left((k_g^2 - \varepsilon k_n^2 + \varepsilon \tau_g^2)^{1/2}\omega\right) - 1}{k_g^2 - \varepsilon k_n^2 + \varepsilon \tau_g^2},$$

and

$$\omega + \int_0^\omega m_{11}ds = \omega + (k_g^2 - \varepsilon k_n^2)\left[\frac{\omega^3}{3!} + (k_g^2 - \varepsilon k_n^2 + \varepsilon \tau_g^2)\frac{\omega^5}{5!} + (k_g^2 - \varepsilon k_n^2 + \varepsilon \tau_g^2)^2 \frac{\omega^7}{7!} - \cdots\right],$$

$$\int_0^\omega m_{12}ds = k_g\left[\frac{\omega^2}{2!} + (k_g^2 - \varepsilon k_n^2 + \varepsilon \tau_g^2)\frac{\omega^4}{4!} + (k_g^2 - \varepsilon k_n^2 + \varepsilon \tau_g^2)^2 \frac{\omega^6}{6!} + \cdots\right]$$
$$- k_n \tau_g \left[\varepsilon \frac{\omega^3}{3!} - \varepsilon(k_g^2 - \varepsilon k_n^2 + \varepsilon \tau_g^2)\frac{\omega^5}{5!} + (k_g^2 - \varepsilon k_n^2 + \varepsilon \tau_g^2)^2 \frac{\omega^7}{7!} + \cdots\right],$$

$$\int_0^\omega m_{13}ds = -\varepsilon k_n\left[\frac{\omega^2}{2!} + (k_g^2 - \varepsilon k_n^2 + \varepsilon \tau_g^2)\frac{\omega^4}{4!} + (k_g^2 - \varepsilon k_n^2 + \varepsilon \tau_g^2)^2 \frac{\omega^6}{6!} + \cdots\right]$$
$$+ \varepsilon k_g \tau_g \left[\frac{\omega^3}{3!} + (k_g^2 - \varepsilon k_n^2 + \varepsilon \tau_g^2)\frac{\omega^5}{5!} + (k_g^2 - \varepsilon k_n^2 + \varepsilon \tau_g^2)^2 \frac{\omega^7}{7!} + \cdots\right].$$

Hence the conditions (i) and (ii) of Theorem 3.2 give that the curve $x(s)$ lying on $S$ fully is periodic of period $\omega$ if and only if $(k_g^2 - \varepsilon k_n^2)^{1/2} = \frac{2k\pi}{i\omega}$, $\tau_g = 0$, where $i$ is imaginary unit with $i^2 = -1$.

## 4. Conclusion

In this paper, by considering the method given by Chung, a criterion for closed space curves in Minkowski space-time $E_v^n$ is given. Furthermore, as a special case, this criterion is applied to a curve lying on a surface in Minkowski 3-space $E_1^3$.